\documentclass[12pt]{amsart}
\usepackage{a4wide}
\usepackage{amssymb,amsthm}
\usepackage{fullpage}
\usepackage{color}
\usepackage{amsmath}

\newcommand{\beq}{\begin{equation}}
\newcommand{\eeq}{\end{equation}}

\Large

\newtheorem{prop}{Proposition}
\newtheorem{theoreme}[prop]{Theorem}
\newtheorem{lem}[prop]{Lemma}

\newtheorem{rem}[prop]{Remark}

\numberwithin{equation}{section}
\numberwithin{prop}{section}
\begin{document}
\title{ A simple criterion of transverse linear instability  \\ for solitary waves}
\author{Frederic Rousset}
\address{IRMAR, Universit\'e de Rennes 1, campus de Beaulieu,  35042  Rennes cedex, France}
\email{ frederic.rousset@univ-rennes1.fr  }
\author{Nikolay Tzvetkov}
\address{D\'epartement de Math\'ematiques, Universit\'e  de Cergy-Pontoise, 95302  Cergy-Pontoise Cedex, France}
\email{nikolay.tzvetkov@u-cergy.fr}
%
%
%
\date{}
\maketitle
\begin{abstract}
We   prove an abstract instability result for an eigenvalue problem with parameter.
We apply this criterion to show the transverse  linear instability of solitary waves
on various examples  from mathematical physics.

\end{abstract}

\section{Introduction}

We shall study  a generalized   eigenvalue problem under the form
\beq
\label{prob}
\sigma A(k) U = L(k) U
\eeq
where  $L(k)$, $A(k)$
are    operators (possibly unbounded)   which depend smoothly on the real parameter  $k$
on some Hilbert space $H$ with moreover  $L(k)$ symmetric.
Our aim is to give an elementary criterion which ensures the existence of $\sigma>0$ and
$k\neq 0$ such that \eqref{prob} has a nontrivial solution $U$. Our motivation for this problem
is the study of transverse instability of solitary waves in Hamiltonian partial differential equations.
 Indeed, let us consider a formally Hamiltonian  PDE, say in $\mathbb{R}^2$,  under the form
  \beq
  \label{pde}
   \partial_{t} \,\mathcal{U} = \mathcal{J} \nabla H( \mathcal{U}), \quad \mathcal{J}^*=-\mathcal{J}
   \eeq
 and  assume that there is a critical point of the Hamiltonian (hence a stationary solution)  $ \mathcal{U}(x,y)= Q(x)$ which depends  only on one variable.  
Note that many equations of mathematical physics have one-dimensional 
 solitary waves solutions which can be seen as critical points of a modified Hamiltonian after
  a suitable change of frame. We shall consider a few examples below. An interesting question
   is  the stability of the one-dimensional state when it is submitted to general two-dimensional
    perturbations.  There are  many examples where  the one-dimensional state  even if it is stable when submitted
     to one-dimensional perturbations is destabilized  by transverse oscillations, we refer for
      example to \cite{Z2}, \cite{APS}, \cite{Kiv}.
    In our previous works \cite{RT1}, \cite{RT2}, \cite{RT3},   we have developed a framework which allows
      to pass from spectral  instability to nonlinear instability. The aim of this note
       is to state a general criterion which allows to get spectral instability.
  Note that  the linearization of \eqref{pde} about $Q$ reads
  \beq
  \label{pdelin}
   \partial_{t} V = \mathcal{J} \mathcal{L} V
   \eeq
   where $\mathcal{L}= D \nabla H(Q)$ is a symmetric operator. Since  $Q$ does not depend on
   the transverse variable $y$, if $\mathcal{J}$ and  $H$ are invariant by translations in $y$,  we can look for a  solution  of \eqref{pdelin}   under the form
   \beq
   \label{fourier} V(t,x,y)= e^{\sigma t } e^{iky } U(x).\eeq
   This yields  an eigenvalue problem for $U$ under the form
   \beq
   \label{probham} \sigma U= (J M) (k) U\eeq
    where   $M(k)$, $J(k)$  defined by
     $$ M(k) U=  e^{-i ky } \mathcal{L}(e^{iky} U), \quad J(k)U =  e^{-i ky } \mathcal{J}(e^{iky} U) $$
    are  operators acting  only in the $x$ variable. 
      Consequently,  if $J(k)$ is invertible, we can set the problem under the form \eqref{prob} with
       $A(k)=  J(k)^{-1}$.   As  we shall see on the examples,
      it may happen that  the skew symmetric operator  $J(k)$ 
      (which very often does not depend on $k$)  is 
      not invertible. 
        In these cases, we can also recast the problem under the form \eqref{prob}. 
         For example, we can look for solutions of \eqref{probham} under the form $U=J(k)^*V$
          and thus get a problem under the form \eqref{prob} with $A(k)=  J(k)^*,$
$L(k) = J(k)M(k) J(k)^*.$

   For the sake of simplicity, we shall work within a real framework but 
   our result can be easily generalized  to   complex  Hilbert  spaces.
     We shall  also study \eqref{prob} only   for $k>0$.
       A similar instability criterion  for $k<0$ can be obtained by  setting  $\tilde{A}(k)= A(-k)$, 
         $\tilde{L}(k)= L(-k)$ and by studying the problem \eqref{prob} for $\tilde{A}$ and $\tilde{L}$.

    Let us fix the functional framework. We consider  a (real) Hilbert space $H$
     with scalar product $(\cdot, \cdot)$. We 
     assume that 
       $L(k)$  is a self-adjoint  unbounded operator with domain  $\mathcal{D}$
        continuously imbedded in  $H$  and independent of
        the real parameter  $k$. Moreover, 
         $L(k)$ as an operator from $\mathcal{D}$ to $H$  is assumed to  depend smoothly on $k$.
        Finally,  we also assume that $A(k)\in \mathcal{L}(\mathcal{D}, H)$ and  depends smoothly
         on $k$. A $\mathcal{C}^1$  dependence  is actually  sufficient for our purpose. 

  Our aim here is to present a criterion which allows to prove transverse instability in solitary
  waves stability problems. This amounts to prove the existence of a nontrivial solution
   of \eqref{prob} with $k \neq 0$ and $\sigma$ with positive real part.    
   In solitary wave stability problem,   $0$ is very often (when the problem is translation invariant in $x$) 
    an eigenvalue of $L(0)$ with eigenvector $Q'$.  Consequently, since we  know
     that  \eqref{prob} has a nontrivial solution for $\sigma=0$, $k=0$,  we can look for a solution
      $(\sigma, U, k)$
      of \eqref{prob} in the vicinity of this particular solution. The main difficulty to implement
       this strategy  is that  also  very often in solitary waves stability problems, 
        $0$ is in the essential spectrum of $JM(0)$, therefore the standard Lyapounov-Schmidt reduction
         cannot be  used. One way to solve this problem is to introduce an Evans function
          with parameter
          $D(\sigma, k)$ (we refer  for example to \cite{AGJ}, \cite{GZ},  \cite{PW}, \cite{KS} for  the definition
           of the Evans function) 
          for the operator $JM(k)$ and then to  study the zeroes of $D$ in the vicinity of $(0,0)$
           (after having proven that  $D$ has  in a suitable  sense a smooth continuation in the vicinity of $(0,0)$).
        We refer for example to \cite{Bridges}, \cite{Benzoni}, \cite{Zumbrun-Serre} 
         for the study of various examples.
          Let us also mention  \cite{chug},   \cite{PS}, \cite{GHS},  for  other  approaches, 
           where   the eigenvalue problem
           is not  reformulated as an ODE with parameters.

        Here  we shall  present a simple approach  which  relies only  on the properties
         of $L(k)$ which are rather easy to check (mostly since it is  a self-adjoint operator) and does not
          rely in principle on the reformulation of the problem as an ODE.

     Our  main assumptions are the following:
     \begin{itemize}
     \item[{\bf (H1)}] There exists $K>0$ and $\alpha>0$ such that $L(k) \geq \alpha\,{\rm Id}$ for $|k| \geq K$;
     \item[{\bf(H2)}]  The essential spectrum  $Sp_{ess}(L(k))$ of $L(k)$
      is included in $[c_{k}, + \infty)$ with $c_{k}>0$ for $k \neq 0$;
\item [\bf{(H3)}]  For every $k_{1} \geq k_{2} \geq 0$,  we have $L(k_{1}) \geq L(k_{2})$. In addition, if for some $k>0$ and $U \neq 0$, we have  $L(k)U= 0$, then
$(L'(k)U, U) > 0$  
(with $L'(k)$  the derivative of $L$ with respect to $k$);
      \item[\bf{(H4)}] The spectrum $Sp(L(0))$ of  $L(0)$ is under the form
        $ \{- \lambda \} \cup I$  where   $- \lambda <0$ is an isolated
         simple eigenvalue   and $I$ is included in $[0, + \infty)$.

     \end{itemize}

          Let us  point out that the structure  of the spectrum of $L(0)$ assumed in (H4) 
          is one of the assumption needed to have the one-dimensional stability of the wave
           (at least when there is a one-dimensional group of invariance in the problem),
           we refer to \cite{GSS}.   Note that  $0$ may be embedded in the essential spectrum
            of $L(0)$.
        \bigskip  

    Our main result is the following:

    \begin{theoreme}
    \label{main}
    Assuming (H1-4), there exists $\sigma >0$, $k \neq 0$ and $U\in \mathcal{D}\backslash\{0\}$
       solutions of \eqref{prob}.
     \end{theoreme}
    Note that we get an instability with $\sigma$ real and positive.
    Once the spectral instability  is established, one may use the general framework developed
     in \cite{RT2} to prove the nonlinear instability of the wave.

The assumption (H3) is clearly matched if   $L'(k)$ is positive  for every  $k>0$.  This last property  is 
  verified for all  the examples that we shall discuss in this paper. Moreover if $L'(k)$ is positive  for $k>0$, the proof of Theorem~\ref{main} can be slightly simplified (see Remark~\ref{referee} below).

    The paper is organized as follows.  
    In the following section, we shall give the proof   of  Theorem \ref{main}.
        Next,  in order to show how our abstract result can be applied, we shall  study  various examples: the KP-I equation, 
       the Euler-Korteweg system and the Gross-Pitaevskii equation. 
        Note that we have already used  similar  arguments  to prove the instability of capillary-gravity
        solitary  water-waves in \cite{RT3}. 
             We hope that  our  approach  can be  useful 
         for other examples,  we also believe that this approach can be adapted  to  many situations
          with slight modifications.

  \section{Proof of Theorem \ref{main}}
      The first step is to  prove that there exists $k_{0}>0$ such that
    $L(k_{0})$ has  a one-dimensional kernel.

     Let us set
     $$ f(k) = \inf_{ \|U \|=1} (L(k) U, U).$$
   Note that   by (H4) $L(0)$ has a negative eigenvalue, hence  
     we have on the one hand  that  
      $f(0)<0$. On the other hand by assumption (H1), we have  that $f(k)>0$ 
      for $k \geq K$.
              Since $f$ is continuous, this implies that there exists a minimal   $k_{0} >0$
       such that $f(k_{0})=0$. 
For every $k<k_{0}$, we get that  $L(k)$ has  a  negative   eigenvalue 
        (since $f(k)$ is negative  and $L(k)$ self-adjoint,  this is a direct consequence of the variational characterization
        of the  bottom of the spectrum    and of (H2) which gives that the essential spectrum of $L (k)$
         is in $(0, + \infty)$). 
        Actually,  there is a unique  negative  simple eigenvalue. Indeed, if we assume
         that $L(k)$ has two (with multiplicity) negative  eigenvalues, then  $L(k)$ is negative 
          on a two-dimensional subspace. By (H3), this yields that  $L(0) \leq L(k)$
           is also negative   on this two-dimensional subspace.
            This contradicts (H4)
            which  contains that  $L(0)$ is nonnegative on  a codimension one subspace.

          By the choice of $k_{0}$ and (H2), we also have that the kernel of $L(k_{0})$ is non-trivial.

         To conclude, we first  note that   
          if for every $k \in (0, k_{0})$ the kernel of $L(k)$ is trivial, then  the kernel of $L(k_{0})$
           is exactly one-dimensional.
           Indeed,  let us pick  $k<k_{0}$, then, since $L(k)$  has a unique simple negative eigenvalue
            and  a trivial kernel, we get that $L(k)$ is positive   on a codimension one  subspace.
             Since $L(k_{0}) \geq L(k)$ by (H3), this implies that  the kernel of $L(k_{0})$
              is exactly one-dimensional. 

                     Next,  we consider the case that  there  exists $k_{1}\in(0, k_{0})$ 
                     such that  $L(k_{1})$ has a nontrivial kernel.
                  Since $L(k_{1})$ has a unique simple negative eigenvalue, we get that 
                  $L(k_{1})$ is nonnegative on a codimension $1$ subspace  $\mathcal{V}= (\varphi)^\perp
                  \equiv \{V\in \mathcal{D}\,:\, (V,\varphi)=0\}$
                  , $\varphi$ being  an  eigenvector associated to the negative eigenvalue.
                   Moreover, thanks to (H2), we  have an orthogonal  decomposition  of $\mathcal{V}$, 
                    \begin{equation}
                    \label{ortho}\mathcal{V}= \mbox{Ker } L(k_{1}) \oplus_{\perp} \mathcal{P}
                    \end{equation}
                     with $\mathcal{P}$ stable
                      for  $L(k_{1})$  and  $L(k_{1})$ restricted to    $\mathcal{P}$  coercive.
                       Note that moreover   $ \mbox{Ker } L(k_{1}) $
                      is  of finite dimension.  For every $U \in \mathcal{S}$ where 
                       $\mathcal{S}$ is the unit sphere of $\mbox{Ker }L(k_{1})$ i.e.
                        $\mathcal{S}=\{ U \in \mbox{Ker }L(k_{1}), \, \|U\|=1\}$, we have  by (H3) 
                      that $ (L'(k_{1})U, U)>0$.
From the compactness of $\mathcal{S}$,  we  get that  $c_{0}= \inf_{U \in \mathcal{S}}
( L'(k_{1})U, U)  $ is positive. This yields that for every $k \geq k_{1}$ close to $k_{1}$
 and $U$ in $\mathcal{S}$,  
 $$
  (L(k) U, U) \geq {c_{0} \over  2} (k- k_{1})$$
   and hence by homogeneity that  
 \begin{equation}
  \label{L'}
  (L(k) U, U) \geq {c_{0} \over  2} (k- k_{1}) \|U\|^2, \quad  \forall U \in \mbox{Ker }L(k_{1}).
 \end{equation}
 Now according to the decomposition \eqref{ortho} of $\mathcal{V}$, we can write
 $ L(k)$  with the block structure
 $$ L(k) = \left(\begin{array}{cc}  L_{1}(k) & A(k) \\ A^*(k) &  L_{2}(k) \end{array}\right).$$
   By the choice of $\mathcal{P}$,  $L_{2}(k_{1})$ is coercive, therefore, there exists $\alpha>0$
    such that for every $k$  close to $k_{1}$, we have
   \begin{equation}
   \label{L2+}( L_{2}(k) U, U) \geq \alpha \|U\|^2, \quad \forall U \in \mathcal{P}.
   \end{equation}
    Moreover,    we also  have that  $A(k_{1})=0$ (since $\mathcal{P}$ is a stable subspace for $L(k_{1})$).
By the assumed regularity with respect to $k$, we thus  get that
\beq
\label{A-}  \|A(k)\|_{\mathcal{L}(\mathcal{P}, Ker \, L(k_{1}) )} \leq M \, |k-k_{1}|, \quad \forall
 k \in [k_{1}/2, 2k_{1}] \end{equation}
  for some $M>0$. 

  Consequently, by using \eqref{L'}, \eqref{L2+} and \eqref{A-}, we get that for every 
  $U=(U_{1},  U_{2}) \in \mathcal{V}$ and every $k>k_{1}$ close to $k_{1}$, we have
  $$ (L(k) U, U) \geq   {c_{0} \over  2} (k- k_{1}) \|U_{1}\|^2+ \alpha \|U_{2}\|^2 -  2 M (k-k_{1}) \|U_{1}\|\, \|U_{2}\|.$$
   From the Young inequality, we  can write
   $$   2 M (k-k_{1}) \|U_{1}\|\, \|U_{2}\| \leq {c_{0} \over 4} \|U_{1}\|^2(k-k_{1}) +  \tilde{M}(k-k_{1}) \|U_{2}\|^2$$
    with $\tilde{M}= 4M^2/c_{0}$ and hence, we obtain
    $$  (L(k) U, U) \geq   {c_{0} \over  4} (k- k_{1}) \|U_{1}\|^2+ ( \alpha - \tilde{M}(k-k_{1})) \|U_{2}\|^2.$$
    In particular, we get that for  every $k>k_{1}$ close to $k_{1}$, $L(k)$ is coercive on $\mathcal{V}$ and
     hence positive. Let us take some $k<k_{0}$ with this last property.
                     Since by (H3), $L(k_{0}) \geq L(k)$, we get
                        that $L(k_{0})$ is also  positive  on  $\mathcal{V}$ which has   codimension $1$.
                         Therefore the kernel of $L(k_{0})$ is exactly one-dimensional.

       We have thus obtained as claimed  that there exists $k_{0}>0$ such that  $L(k_{0})$ has a one-dimensional
        kernel. Thanks to (H2), we  also have that $L(k_{0})$ is a  Fredholm  operator with zero index.
        We can therefore use the Lyapounov-Schmidt method to study the eigenvalue problem
         \eqref{prob} in the vicinity of $\sigma=0$, $k=k_{0}$ and $U= \varphi$ where
          $\varphi$ is in the kernel of $L(k_{0})$ and  such that $\| \varphi \|=1$.

   We look for $U$ under the form  $U=\varphi+V$, where 
$$
V\in {\varphi}^{\perp}\equiv \{V\in \mathcal{D}\,:\, (V,\varphi)=0\}.
$$
Therefore we need to solve $G(V,k,\sigma)=0$ with $\sigma>0$, where
$$
G(V,k,\sigma)=L(k)\varphi+L(k)V-\sigma A(k) \varphi-\sigma A(k) V,\quad V\in {\varphi}^{\perp}\,.
$$
We shall use the implicit function theorem to look for $V$ and $k$ as functions of $\sigma$.
Note that the same approach is for example used in \cite{GHS}.
We have that
\beq
\label{DG}
D_{V,k}G(0,k_0,0)[w,\mu]=\mu\Big[\frac{d}{dk}L(k)\Big]_{k=k_0}\varphi+L(k_0)w\,.
\eeq
By using (H3),  we obtain that $D_{V,k}G(0,k_0,0)$ is a
bijection from $ {\varphi}^{\perp} \times \mathbb{R}$ to $H$. 
We can thus  apply the implicit function theorem to get that
for $\sigma$ in a neighborhood  of  zero there exists $k(\sigma)$ and $V(\sigma)$ such that
$ G(V(\sigma), k(\sigma), \sigma)=0$.  This ends
  the proof of Theorem~\ref{main}.  
\begin{rem}\label{referee}
Let us remark that if we assume  that $L'(k)$ is positive   for $k>0$ in place of (H3),  then we can simplify  the argument 
giving a $k_0\neq 0$ such that $L(k_0)$ has a one-dimensional kernel.
Namely, in this case by using (H4), we  have
           that $L(0)$ is nonnegative on a codimension $1$ subspace $\mathcal{V}$
           (given by $\mathcal{V}= \pi_{[0, + \infty)}(L(0)) \cap \mathcal{D}$, where  $\pi_{[0, + \infty)}(L(0))$
            is the spectral projection on the nonnegative spectrum of $L(0)$).
Next, using that $L'(s)$ is positive   for $s>0$, we get for every $k>0$  that
             $$ (L(k)U, U) = \int_{0}^k (L'(s)U, U)\, ds + (L(0)U, U) \geq 0, \quad \forall U \in \mathcal{V}.$$
               Moreover,  if   $(L(k)U, U)= 0$ for $U \in \mathcal{V}$ then the above identity yields
             $$   \int_{0}^k (L'(s)U, U)\, ds=0$$
              and hence  by using again that  $L'(k)$ is positive  for $k>0$, 
we obtain that $U=0$.
               Consequently, we get that  for $k>0$,  $L(k)$ is positive on a codimension $1$ subspace.
            This yields that the dimension of the kernel of $L(k_{0})$ is exactly one.
             \end{rem}

\section{Examples}

In this section we shall  study various physical examples where Theorems \ref{main} 
can be used to prove  the instability of line solitary waves.
\subsection{KP-I equation}
We shall first see that the  instability argument given in \cite{RT2} can be interpreted
in the framework of \eqref{prob}.
Let us consider the generalized KP-I equation where
\beq
\label{KP} \partial_{t} u = \partial_{x}\big( - \partial_{xx} u  - u^p \big) + \partial_{x}^{-1}
 \partial_{yy} u , \quad p=2,\, 3, \, 4
\eeq
and $u(t,x,y)$ is real valued. There is an explicit one-dimensional solitary wave
 (which thus  solves the generalized  KdV equation):
  $$ u(t,x)= Q(x-t) = \Big({ p+ 1 \over 2} \Big)^{1 \over p-1}  \Big(
   \mbox{sech}^2 \big( {(p-1) (x-t) \over 2 } \big) \Big)^{ 1 \over p-1 }.$$
Note that in this problem, it suffices to study the stability of the speed one solitary wave
since the solitary wave with speed $c>0$ can be deduced from it by scaling:  the solitary wave
 with speed $c>0$ is given by
 $$ Q_{c}(\xi)=  c^{ 1 \over p-1} Q(\sqrt{c}\, \xi).$$

After   changing $x$ into  $x-t$ (and keeping the notation $x$) and linearizing about $Q$,
we can seek for solution under the form 
$$ e^{\sigma t } e^{iky } V(x)$$
to get
the equation
$$\sigma  V= \partial_{x} \big( - \partial_{xx}  - k^2 \partial_{x}^{-2}+  1 - p Q^{p-1} \big) V.$$
 We can seek for a solution $V$ under the form $V= \partial_{x} U$
  to get that $U$ solves
$$  - \sigma  \partial_{x} U=\Big(-  \partial_{x}(  - \partial_{xx}  +1 -   p Q^{p-1} ) \partial_{x}  +k^2 \Big)U.$$
 Therefore, this eigenvalue problem is under the form \eqref{prob} with
 $$ A(k)= - \partial_{x}, \quad L(k)
  =-  \partial_{x}( -  \partial_{xx}  + 1 -  p Q^{p-1} ) \partial_{x}  +k^2.$$
    By  choosing $H= L^2(\mathbb{R})$ and $\mathcal{D}= H^4(\mathbb{R})$, we are in an appropriate functional framework.
  Note that $L(k)$ has a self-adjoint realization.

 Let us check the assumptions (H1-4).

  Since we have  
  $$ (L(k) U, U) \geq   \| \partial_{xx} U \|_{L^2}^2 + k^2 \|U \|_{L^2}^2  + \| \partial_{x} U \|_{L^2}^2
   -  p \|Q^{p-1} \|_{L^\infty } \| \partial_{x} U \|_{L^2}^2$$
    and that
   for every $\delta >0$, there exists $C(\delta) >0$ such that 
  $$  \| \partial_{x} U \|_{L^2}^2 \leq \delta \| \partial_{xx} U \|_{L^2}^2 + C(\delta) \|  U \|_{L^2}^2,$$
   we immediately get that (H1) is verified.

  Next, we note that  $L(k)$ is a compact perturbation of 
  $$L_{\infty}(k)=   -  \partial_{x}( -  \partial_{xx}  + 1  ) \partial_{x}  +k^2,  $$ 
  we  thus get from the Weyl Lemma and  the explicit knowledge of the spectrum of $L_{\infty}(k)$ that
   the  essential spectrum of $L(k)$ is included in $[k^2, + \infty)$ and thus 
    that (H2) is verified.

    Assumption (H3) is obviously verified since $L'(k)$ is positive for every $k>0$.

 Finally, let us check (H4). Note that $L(0)= -\partial_{x} C \partial_{x}$ where
  $C$ is  a second order differential operator. We notice that
   $C  Q' =0$ and that by the same argument as above, the essential spectrum of $C$
    is contained in $[1, + \infty)$. Since $ Q'$ vanishes only once, we get by
     Sturm Liouville theory (we refer for example to \cite{Dunford}, chapter XIII) that $C$ has  a unique negative eigenvalue with associated eigenvector $\psi$.
      Moreover, we also have that
     \beq
     \label{Calpha}
      (Cu, u ) \geq  0
       \quad \forall u \in (\psi)^{\perp}
      \eeq
      After these preliminary remarks, we can  get  that
       $L(0)$ has   a  negative eigenvalue. Indeed
        by an approximation argument, we can construct a sequence $u_{n}$ in $\mathcal{D}$
         such that $\partial_{x} u_{n}$ tends to  $\psi $ in $\mathcal{D}$ then, for
          $n$ sufficiently large $(L(0) u_{n}, u_{n}) = (C\partial_{x} u_{n}, \partial_{x}u_{n})$ is negative.
        By the variational characterization of the lowest eigenvalue, we get that
         $L(0)$ has a negative eigenvalue.
          Moreover, for every  $U$ such that $ (\partial_{x}U, \psi)=0$, 
           we have that  
          $$ \big(L(0) U, U\big) = \big(C \partial_{x} U, \partial_{x} U\big)\geq0$$
       This proves that $L(0)$ is nonnegative  on a codimension one  subspace
        and hence that there is at most one negative eigenvalue.
        We have thus proven that (H4) is verified.

         Consequently, we get from Theorem \ref{main} that the  solitary wave 
         is transversally unstable.

\subsection{Euler-Korteweg models}
We consider  a general class of  models describing the isothermal motion of compressible fluids
and taking into account  internal capillarity.  
The main feature of these models is that the free energy $F$ depends both on
 $\rho$ and $\nabla \rho$. In the isentropic case,  we have: 
 $$ F(\rho, \nabla \rho) = F_{0}(\rho) + {1 \over 2} K(\rho)  |\nabla \rho |^2$$
  where $F_{0}(\rho)$ is the standard part and $K(\rho)>0$ is a capillarity coefficient.  The pressure
   which is defined by 
$p= \rho {\partial F \over \partial \rho} - F$ 
reads     
$$ p(\rho, \nabla \rho)=  p_{0}(\rho) + {1 \over 2} \big( \rho K'(\rho) - K(\rho)  \big) | \nabla \rho|^2  \big.
$$
The equations of motion read
 \begin{eqnarray}
 \label{euler1}
& &   \partial_{t} \rho + \nabla \cdot( \rho u ) = 0 , \\
& & \label{euler2}   \partial_{t} u + u \cdot \nabla u + \nabla (g_{0}(\rho) ) =  \nabla \big( K(\rho) \Delta \rho +
  {1 \over 2} K'(\rho ) | \nabla \rho |^2 \big).
  \end{eqnarray}
  In this model, $\rho>0$ is the density of the fluid and $u$ the velocity, 
  $g_{0}$ (which is linked to $p_{0}$ by 
   $\rho g_{0}'(\rho)=p_{0}'(\rho)$) 
   and $K(\rho)>0$ are smooth functions of $\rho$ for $\rho>0$.

We shall consider a one-dimensional solitary wave  of \eqref{euler1}, \eqref{euler2} under the form
$$(\rho(t,x,y), u(t,x,y))= (\rho_{c}(x-ct), u_{c}(x-ct))= Q_{c}(x-ct)$$
 such that
 \beq
 \label{rhoinfty}
  \lim_{x \rightarrow \pm \infty} Q_{c}=  Q_{\infty}=(\rho_{\infty}, u_{\infty}), \quad \rho_{\infty}>0.
  \eeq
     We shall assume that
        \beq
    \label{hypeuler}
     \rho_{\infty} g_{0}'(\rho_{\infty}) > (u_{\infty}- c)^2.
     \eeq
     This condition ensures  that  $Q_{\infty}$ is a saddle  point in the ordinary differential equations
      satisfied by the profile. Under this condition,  one can find solitary waves, moreover, 
       they have the interesting property   that
    $ \rho_{c}'$ vanishes only once. We refer for example to \cite{Benzoni-Danchin},
     for the study  of the existence of
    solitary waves for this system.

    Here we shall study the (linear)  transverse instability of these solitary waves.
     We shall restrict our study to  potential solutions of  \eqref{euler1}, \eqref{euler2}
      that is to say solutions  such that $u= \nabla \varphi$. Note that this will give a better
       instability result, this means that we are able to find instabilities even in the framework of potential
        solutions. 

      This yields the system 
      \begin{eqnarray}
 \label{euler3}
& &   \partial_{t} \rho + \nabla\varphi \cdot \nabla \rho + \rho \Delta \varphi = 0 , \\
& & \label{euler4}   \partial_{t} \varphi + {1 \over 2} | \nabla \varphi|^2 +  g_{0}(\rho)  =  K(\rho) \Delta \rho +
  {1 \over 2} K'(\rho ) | \nabla \rho |^2.
  \end{eqnarray} 
 Changing $x$ into $x-ct$ (and keeping the notation $x$) to make the wave stationary,  linearizing
  \eqref{euler3}, \eqref{euler4} about a solitary wave $Q_{c}= (\rho_{c}, u_{c})$ and
   looking for solutions $(\eta, \varphi)$ under the form
   $$ (\eta, \varphi)= e^{\sigma t } e^{iky } U(x), $$
   we find an eigenvalue problem
    under the form \eqref{prob} with $A(k)= J^{-1}$ and 
    $$ J= \left( \begin{array}{cc} 0 & 1 \\ -1 & 0 \end{array} \right),$$
    $$ L(k)= \left( \begin{array}{cc}
      - \partial_{x}\big( K(\rho_{c}) \partial_{x})  +k^2 K(\rho_{c})  -  m &  -c \partial_{x} + u_{c} \partial_{x} \\
        c \partial_{x} - \partial_{x}\big( u_{c}\cdot \big) & - \partial_{x}\big(\rho_{c} \partial_{x} \cdot \big)
         + \rho_{c}k^2 \end{array} \right)$$
         where the function $m(x)$ is defined by
         $$ m= K'(\rho_{c}) \rho_{c}'' + {1 \over 2} K''(\rho_{c}) (\rho_{c}')^2 - g_{0}'(\rho_{c}).$$
 By taking $H= L^2(\mathbb{R}) \times L^2(\mathbb{R})$, $\mathcal{D}= H^2 \times H^2$, 
  we are in the right functional framework, in particular,   $L(k)$ has a self-adjoint realization. Let us now check assumptions
    (H1-4):
    \begin{itemize}
    \item (H1):  with $U= (\rho, \varphi)$, we have
 \begin{multline*}
  \big(L(k)U, U \big) \geq  \int_{\mathbb{R}} \Big( K(\rho_{c}) \big( | \partial_{x} \rho |^2 +
k^2 |\rho|^2  \big)  + \rho_{c} \big(| \partial_{x} \varphi |^2 + k^2
 | \varphi |^2  \big) \\
  - \mathcal{O}(1)\big( |\rho |( |\rho | + | \varphi |)  +  | \partial_{x} \rho | \, | \varphi| \big)\Big)\,dx
   \end{multline*}
  where $\mathcal{O}(1)$ is independent of $k$.      
    Since $K(\rho_{c}) \geq \alpha >0$, we get 
   by using the Young inequality:
      \beq
      \label{young}
     ab \leq {\delta \over 2} a^2 + { 1 \over 2  \delta } b^2, \quad a,\, b \geq 0, \quad \delta >0,\eeq
            that (H1) is verified for $k$ sufficiently large.
 \bigskip           

         \item (H2)    By standard arguments (see \cite{Henry}, for example), 
          to locate the essential spectrum  of $L(k)$,  we have  to study  the 
   spectrum  of
$$   L_{\infty}(k) =     
\left( \begin{array}{cc}
       \big( K(\rho_{\infty})( -  \partial_{xx}  +k^2 \big)  + 
        g_{0}'(\rho_{\infty})&  -c \partial_{x} + u_{\infty} \partial_{x} \\
        c \partial_{x} - u_{\infty } \partial_{x} &  \rho_{\infty} \big( - \partial_{xx} + k^2\big)
        \end{array}
          \right).$$
          By using the Fourier transform, we  can compute explicitly the spectrum of this operator, 
           we find that $\mu$ is in the spectrum  of $L_{\infty}(k)$ if  and only if there exists $\xi$
            such that
        $$ \mu^2 - s(\xi, k) \mu + p(\xi, k)= 0$$
         with
        \begin{eqnarray*}
       & & s=  (K(\rho_{\infty})  + \rho_{\infty})(k^2+ \xi^2) + g_{0}'(\rho_{\infty} ), \\
        & & p= 
            \rho_{\infty} K(\rho_{\infty})(k^2 + \xi^2 )^2  + \rho_{\infty} g_{0}'(\rho_{\infty}) k^2 
           +  \big(\rho_{\infty }  g_{0}'(\rho_{\infty})  - (u_{\infty} - c )^2 \big) \xi^2 \geq 0 .\end{eqnarray*}
          By using that $\rho_{\infty}$ and $K(\rho_{\infty})$ are positive and the condition \eqref{hypeuler},
           we get  that the two roots are nonnegative for all $k$ and  strictly positive for $k \neq 0$.
          This proves that  (H2) is matched.
      \bigskip

   \item(H3)  We have
   $$L'(k) = \left( \begin{array}{cc}   2 k K(\rho_{c})  & 0    \\  0 & 2 \rho_{c} k \end{array} \right).$$
    Consequently, (H3)  is verified since $\rho_{c}$ and $K(\rho_{c})$ are positive.
    \bigskip

    \item (H4) We shall use the following algebraic lemma:
    \begin{lem}
    \label{lemalg}
Consider a symmetric operator  on $H$ under the form
$$ L= \left( \begin{array}{cc} L_{1} & A \\ A^*  & L_{2} \end{array}\right)$$
 with   $L_{2}$  invertible.  Then we have
 $$ (LU, U)= \Big( \big( L_{1}-  A L_{2}^{-1} A^*\big) U_{1}, U_{1} \Big) +  \Big( L_{2}
  \big( U_{2}+ L_{2}^{-1} A^* U_{1} \big), U_{2} + L_{2}^{-1} A^*U_{1} \Big).$$
     \end{lem}
     The proof is a direct calculation.  Note that 
       the above lemma remains true as soon as the quadratic form in the right-hand side 
        makes sense  (and hence even if $L_{2}^{-1}$ is not well-defined.)

     Let us apply this lemma to    $L(0)$. We see that with
    $$ A=  (u_{c} - c ) \partial_{x}, \quad L_{2}= - \partial_{x} (\rho_{c} \partial_{x} \cdot),$$
     if $u\in H^2$  solves the equation
     $L_{2}  u = -  A^* U_{1}$,   then  
     $$   \partial_{x} u=  -  {1 \over \rho_{c} } (u_{c} - c) U_{1}.$$
Consequently, we get
$$  AL_{2}^{-1} A^* U_{1} = (u_{c}- c ) \partial_{x} u =   {(u_{c} -c )^2 \over \rho_{c}} U_{1}$$
and hence we have the following factorization:
\beq
\label{factor1} (L(0) U, U)= (MU_{1}, U_{1}) + 
\int_{\mathbb{R}} \rho_{c}  \Big| \partial_{x} U_{2}  + {1 \over \rho_{c}} (u_{c} - c) U_{1} \Big|^2\, dx \eeq
 where 
$$  MU_{1}=    -\partial_{x}\big( K(\rho_{c}) \partial_{x} U_{1}\big)  -  m \,U_{1} - {(u_{c} - c)^2 \over \rho_{c}}
 U_{1}.$$
 Note that $M$ is a second order differential operator and that by using  the profile equation satisfied
  by $Q_{c}$, we can check that  $\rho_{c}'$ is in the kernel of $M$.  Since $\rho_{c}'$
   has a unique zero,  this proves that $M$ has exactly one negative eigenvalue with corresponding
    eigenfunction $R$.
   From  the condition \eqref{hypeuler}, we also get that the essential spectrum of $M$
    is included in $[ \alpha, +\infty)$ for some $\alpha>0$. In particular (since $M$ is self-adjoint), we get that
    \beq
    \label{M1}
     (MU_{1}, U_{1}) \geq  0,  \quad \forall \,U_{1} \in (R)^\perp.
     \eeq

We can now use these properties of $M$ to prove that (H4) is matched.
     We can first get  from \eqref{factor1} that $L(0)$ has indeed one negative direction. A first try
      would be to take $U_{1}= R$ and 
      $$ \partial_{x} U_{2}= -  {1 \over \rho_{c}} (u_{c} - c) R.$$
     The problem is that this equation does not have a solution in  $L^2$.  Nevertheless, we can get the
      result by using an approximation argument. Indeed, again by cutting the low frequencies,  we can  choose a sequence  $U_{2}^n  \in H^2$
       such that 
        $$\partial_{x} U_{2}^n  \rightarrow  -   {1 \over \rho_{c}} (u_{c} - c)R$$
        in $H^2$.
        Then since $(MR, R)<0$, we get that
       $  (L(0) U^n, U^n)<0$, with $U^n= (R, U_{2}^n)$, for $n$ sufficiently large
        and hence by the variational characterization of the smallest eigenvalue, we get that
         $L(0)$ has a negative eigenvalue. 
          From \eqref{factor1} and \eqref{M1}, we then get that this negative eigenvalue is unique.
         This proves that (H4) is verified. 
       \end {itemize}

  Consequently, we can use  Theorem~\ref{main}
    to get the instability of the solitary wave.
   We have thus proven:
   \begin{theoreme}
   \label{theoEK}
  If a solitary wave satisfies the condition \eqref{hypeuler}, then it is unstable
  with respect to transverse perturbations.

   \end{theoreme}
   Note that a similar result has been obtained in \cite{Benzoni} by using an Evans function 
    calculation.

\subsection{Travelling waves  of the Gross-Pitaevskii equation}

In this subsection, we consider the Gross-Pitaevskii equation which is  a standard model
for Bose-Einstein condensates, 
\beq
\label{GP}
i \partial_{t} \psi + { 1 \over 2 }\Delta \psi + \psi (1 - | \psi|^2) =0 \eeq
where  the unknown $\psi$ is complex-valued. This equation has well-known explicit
one-dimensional   travelling waves (the so-called dark solitary waves) whose modulus tend to $1$
at infinity, for every $c < 1 $ they read :
\beq
\label{dark} \psi(t,x,y)= \Psi_{c}(z) = 
 \sqrt{ 1 - c^2 }\, \mbox{tanh}\,\big(  z  \sqrt{ 1 - c^2} \big) \big) + i c , \quad z= x-ct.\eeq
In the case of the standard solitary waves of the cubic focusing Schr\"odinger equation,
their transverse instability which was shown by
 Zakharov and Rubenchik  can be studied by a  standard bifurcation argument since
 $0$ is not in the  essential  spectrum of the linearized operator, we refer  for example to 
  \cite{RT1} for the details. This is  not the case  for the dark solitary waves, $0$ is in the essential spectrum of the linearized 
   operator, we shall thus use the criterion given by Theorem \ref{main}.   

Note that for $c \neq 0$,  $\Psi_{c}$ does not vanish. Consequently, we can study
 the stability of these waves (travelling bubbles) by using the  Madelung transform, i.e. by seeking solutions
  of  \eqref{GP} under the form
  $$ \psi = \sqrt{\rho}e^{i \varphi}$$
  with smooth  $\rho$ and $\varphi$. We then classically find that
   $(\rho, u= \nabla \varphi)$ is  a solution of   \eqref{euler1}, \eqref{euler2}  with: 
   $$ g_{0}(\rho)= \rho-1, \quad K(\rho) = {1 \over 4 \rho}.$$
   The dark solitary waves  for $c \neq 0$ becomes a solitary wave $(\rho_{c} , u_{c})$
    of \eqref{euler1}, \eqref{euler2}  with
    $$ \rho_{c}(z)=  c^2 + (1-c^2) \mbox{tanh}^2 \,\big(  z  \sqrt{ 1 - c^2} \big)\Big), \quad$$
    $$ u_{c}(z)= - {c (1 -c^2) \over \rho_{c} }   \Big(  1 - \mbox{tanh}^2 \,\big( { z } \sqrt{ 1 - c^2} \big) 
     \Big).$$
    In particular, we thus have $\rho_{\infty}=  1$ and $u_{\infty}= 0$.  Since $g_{0}'(\rho)=1$, the condition
     \eqref{hypeuler} reduces to $c^2<1$. Consequently,  we have by Theorem \ref{theoEK} 
       that all the dark solitary waves with
       $ |c|<1,$ $c \neq 0$ are unstable with respect to transverse perturbation.
        Note that the one-dimensional stability  of these travelling bubbles was  shown  in \cite{Lin}. 

     It remains to study the case $c=0$. Note that $\Psi_{0}$ is a stationary solution, the so-called  black soliton, which has
      the very simple expression
      $$ \Psi_{0}(x)= \mbox{tanh }\big(x).$$
      Its one dimensional orbital stability has been shown in \cite{Gerard},  \cite{Gravejat-Saut}.
       Here we shall prove that it  becomes  transversally unstable by using Theorem \ref{main}.
        Since the Madelung transform is not appropriate (the solitary wave vanishes
         at the origin)  we shall work directly on 
         the formulation \eqref{GP}.

        Linearizing \eqref{GP} about $\Psi_{0}$,  splitting real and imaginary parts and
         seeking solutions under the form  \eqref{fourier} yield a problem
         under the form \eqref{prob} with $A(k)= J^{-1}$ and
      $$ J= \left( \begin{array}{cc} 0 & 1 \\ -1  & 0 \end{array}\right), \quad  
  L(k) = \left( \begin{array}{cc}     {1 \over 2 }( -  \partial_{xx} + k^2) + 3 \Psi_{0}^2 - 1 & 0 \\ 0
       &  {1\over2}(-\partial_{xx} + k^2)  -  (1-\Psi_{0}^2)  \end{array}\right).$$
 Again, with $H= L^2 \times L^2$ and $\mathcal{D}= H^2 \times H^2$, we can
  check (H1-4).

  (H1)  and (H3)  are obviously matched. Thanks to the decay of the solitary wave, we
  have that  $L(k)$ is a compact perturbation of 
  $$ L_{\infty}(k)= \left( \begin{array}{cc}     {1 \over 2 }( -  \partial_{xx} + k^2) +  2  & 0 \\ 0
       &  {1\over2}(-\partial_{xx} + k^2)    \end{array}\right)$$
and hence, a simple computation shows that (H2) is also matched.
Finally,  we can also easily   check (H4).
   Let us set $L(0)= \mbox{diag }(L_{1}, L_{2})$. We first notice
that 
$$ L_{1}  \Psi_{0}'=0, \quad L_{2} \Psi_{0}=0.$$
and that  the essential spectrum of $L_{1}$ is contained in $[2, + \infty)$ and the one
 of $L_{2}$ in $[0, +\infty)$.
  Since $\Psi_{0}'$ does not vanish and $\Psi_{0}$ vanishes only once, we get by Sturm-Liouville theory that 
   $0$ is the first eigenvalue of $L_{1}$ and  that 
         $L_{2}$ has a unique negative eigenvalue. 
    This proves that (H4) is  matched. 

      Consequently, we get from Theorem \ref{main} that the black soliton  $\Psi_{0}$ is transversally unstable.  

    We have thus proven:
    \begin{theoreme}
    For every $c$, $|c|<1$, the dark solitary waves \eqref{dark} are transversally unstable.

    \end{theoreme}

\begin{rem}
Using arguments as above, we can also  prove the transverse instability of the one dimensional localized solitary waves of the nonlinear Schr\"odinger equation and thus 
obtain another proof of the classical  Zakharov-Rubenchik instability result.
\end{rem}

\begin{rem}
The most difficult assumption to check is often the assumption (H4). Note that on the above
 examples this is always a direct  consequence of  Sturm-Liouville theory  which is an ODE result.
In the above examples, the eigenvalue problem for $L(0)$ is  already itself an ODE.
 Nevertheless,   for the capillary-gravity solitary waves problem  studied in \cite{RT3}, 
there  is a  nonlocal operator arising  in the definition of  $L(0)$ 
   and hence the eigenvalue  problem  for $L(0)$  cannot be formulated as an ODE.
   Nevertheless, it was  proven  by Mielke \cite{Mielke} that
    (H4) is matched since  in the KdV limit the spectral properties of $L(0)$ are the
    same  as the ones of the linearized KdV hamiltonian about the KdV solitary wave
     which are known (again thanks to Sturm Liouville theory).
\end{rem}

  \vspace{1cm}
  {\bf Acknowledgements.}
   We thank Jean-Claude Saut and David Chiron for fruitful discussions about  this work. 
   We also warmly thank the referee  for  the careful reading of the manuscript and many remarks
     which
   have  greatly improved the result and the presentation.

\end{document}